\documentclass[12pt]{article}
\usepackage[leqno]{amsmath}
\usepackage{amsfonts}
\usepackage{graphicx}
\usepackage{amsmath}
\usepackage{amssymb}
\usepackage{latexsym}
\usepackage{amsmath,amsfonts,amssymb,amsthm,euscript,makeidx,color,mathrsfs}
\usepackage{boxedminipage}
\usepackage{xcolor}

\def\sqr#1#2{{\vcenter{\vbox{\hrule height.#2pt
              \hbox{\vrule width.#2pt height#1pt \kern#1pt \vrule width.#2pt}
              \hrule height.#2pt}}}}

\def\3n{\negthinspace \negthinspace \negthinspace }
\def\2n{\negthinspace \negthinspace }
\def\1n{\negthinspace }

\def\={\buildrel \triangle \over =}

\def\D{\Delta}

\def\cA{{\cal A}}

\def\cQ{{\cal Q}}

\def\cS{{\cal S}}

\def\ms{\medskip}

\def\deq{\mathop{\buildrel\D\over=}}

\def\({\Big (}
\def\){\Big )}
\def\[{\Big[}
\def\]{\Big]}

\def\bde{\begin{definition}\label}
\def\ede{\end{definition}}
\def\be{\begin{equation}}
\def\bel{\begin{equation}\label}
\def\ee{\end{equation}}
\def\bt{\begin{theorem}\label}
\def\et{\end{theorem}}
\def\bc{\begin{corollary}\label}
\def\ec{\end{corollary}}
\def\bl{\begin{lemma}\label}
\def\el{\end{lemma}}
\def\bp{\begin{proposition}\label}
\def\ep{\end{proposition}}
\def\bas{\begin{assumption}}
\def\eas{\end{assumption}}
\def\br{\begin{remark}\label}
\def\er{\end{remark}}
\def\bex{\begin{example}\label}
\def\ex{\end{example}}

\newtheorem{theorem}{Theorem}[section]
\newtheorem{corollary}[theorem]{Corollary}

\newtheorem{lemma}[theorem]{Lemma}
\newtheorem{proposition}[theorem]{Proposition}

\theoremstyle{definition}
\newtheorem{definition}[theorem]{Definition}
\newtheorem{example}[theorem]{Example}
\newtheorem{remark}[theorem]{Remark}

\def\ba{\begin{array}}
\def\ea{\end{array}}
\def\ben{\begin{enumerate}}
\def\een{\end{enumerate}}
\def\ed{\end{document}}

\def\ges{\geqslant}

\def\square#1{\vbox{\hrule\hbox{\vrule height#1%
     \kern#1\vrule}\hrule}}
\def\rectangle#1#2{\vbox{\hrule\hbox{\vrule height#1%
     \kern#2\vrule}\hrule}}

\font\tenbb=msbm10 \font\sevenbb=msbm7 \font\fivebb=msbm5

\newfam\bbfam
\scriptscriptfont\bbfam=\fivebb \textfont\bbfam=\tenbb
\scriptfont\bbfam=\sevenbb

\newtheorem{assumption}{Assumption}[section]

\makeatletter
   
   \@addtoreset{equation}{section}
\makeatother

\begin{document}
\title{Numerical Method for a Class of Algebraic Riccati Equations}

\author{Lucky Qiaofeng Li\footnote{Li Po Chun United World College, Hong Kong. Email: luckyql17@gmail.com.}, \quad
Xizhi Su\footnote{Department of Mathematics, The National University of Singapore, Republic of Singapore. Email: suxizhi@u.nus.edu.}}

\date{}

\maketitle

\begin{abstract}
We study an iteration approach to solve the coupled algebraic Riccati equations when they appear in general two player closed-loop type Nash differential games over an infinite time horizon.
Also, we propose an effective algorithm for finding positive definite solutions.
In particular, we present various numerical examples connected with matrix Riccati equations according to different dimensions.
\end{abstract}


%

\section{Introduction}

Before we present our research, it is worth pointing out the real life experience (community volunteer activity) relevant to our research work.
Besides the interest in our own problem, this class of equations appears, for instance, in the optimal control problem for infinite Markov jump linear systems \cite{BF} or linear-quadratic two-person zero-sum differential games over an infinite horizon \cite{SYZ}.
An important recent development in the use of invariant subspace-based methods is the use of iterative refinement. As we know, this idea is closely related to the solution of the Riccati equation in Newton's method. In particular, Newton's method itself is an interesting yet useful technique, and now a fairly complete theory has been developed for this algorithm, including a computable estimate of the convergence area. We refer the interested reader to \cite{Laub,GH,GL,BS,DP} for details.

The algebraic Riccati equations play fundamental roles in engineering, management, economic, finance, linear-quadratic control and estimation systems as well as in many branches of applied mathematics. The purpose of this article is not to investigate the vast literature available for these equations, but to provide readers with references \cite{Laub,GH,GL,BS,DP,XL,KS,AJ,LSY}.
Generally, algebraic Riccati equations that appear in game theory are non-symmetric, so this makes the effective numerical methods developed for symmetric equations based on control theory inapplicable. Our problem is to find the solution of the coupled algebraic Riccati equations generated by two-person linear-quadratic non-zero sum deterministic differential games over an infinite horizon with closed-loop Nash equilibria. These equations based on closed-loop type Nash equilibria are symmetric, but they are coupled together so that the methods established to solve control problems cannot be directly applied.
In this article, we develop an iteration approach to tackle these coupled algebraic Riccati equations, and propose an effective algorithm for finding positive definite solutions.

The article is organized as follows. Section 2 introduces some basic notions and presents formulation and algorithm.  In Section 3, we present several numerical examples according to different dimensions to illustrate the proposed approach. Section 4 concludes this research work.

\section{Formulation and Algorithm}\label{sect:algorithm}

Throughout this paper, we let $\mathbb{R}^{n\times m}$
and $\mathbb{S}^n$ be the set of all $(n\times m)$ (real) matrices and $(n\times n)$ symmetric (real) matrices.
To simplify the notation, 
we denote (for $i=1,2$)
$$
\left\{\begin{array}{ll}
A \in \mathbb{R}^{n\times n}, \quad 
B = (B_1, B_2), \quad B_1 \in \mathbb{R}^{n\times m_1}, \quad B_2 \in \mathbb{R}^{n\times m_2}, \\ [2mm]
Q_1, Q_2 \in \mathbb{S}^n, \quad 
R_i=\begin{pmatrix} R_{i11}&R_{i12}\\ R_{i21}& R_{i22}\end{pmatrix}, \\ [5mm] 
R_{i11} \in \mathbb{S}^{m_1},  \quad R_{i22} \in \mathbb{S}^{m_2},  \quad R_{i12}=R_{i21}^\top \in \mathbb{R}^{m_1\times m_2},
\end{array}\right.
$$
where the notation $M^\top$ denotes transpose of matrix $M$.

\begin{definition} \it
The following is called a system of coupled {\it algebraic Riccati equations} (AREs):
\begin{equation}\label{eq:Riccati}
\left\{\begin{array}{l}
X_1A+A^\top X_1+\Theta^\top R_1\Theta + X_1B\Theta + \Theta^\top B^\top X_1 + Q_1 = 0, \\ [2mm]
X_2A+A^\top X_2+\Theta^\top R_2\Theta + X_2B\Theta + \Theta^\top B^\top X_2 + Q_2 = 0,
\end{array}\right.
\end{equation}
where $\Theta = (\Theta_1, \Theta_2)^\top$, $\Theta_1 \in \mathbb{R}^{m_1\times n}$ and $\Theta_2 \in \mathbb{R}^{m_2\times n}$ satisfies
\begin{equation}\label{eq:cond}
\begin{array}{l}
\begin{pmatrix} B_1^\top X_1 \\ B_2^\top X_2\end{pmatrix}
+ \begin{pmatrix} R_{111} & R_{112} \\ R_{221} & R_{222}\end{pmatrix} 
\begin{pmatrix}\Theta_1 \\ \Theta_2\end{pmatrix} = 0, 
\end{array}
\end{equation}
and
$$
R_{111} \succ 0, \quad R_{222} \succ 0.
$$
\end{definition}

\ms

\ms

\textbf{Problem.} 
For $i=1,2$, $Q_i \succeq 0$ and $R_i \succ 0$, find a pair of solution $(X_1, X_2)$ to satisfy equations (\ref{eq:Riccati})-(\ref{eq:cond}), where $X_i \in \mathbb{S}^n$ and $X_i \succeq 0$.

\ms

\ms

Since we cannot solve the coupled algebraic Riccati equations (\ref{eq:Riccati})-(\ref{eq:cond}) directly, we separate the coupled equations into two decoupled equations (\ref{ARE-1}) and (\ref{ARE-2}) for finding positive definite solutions using iterations in the following algorithm.

\ms

\ms

\begin{boxedminipage}[t]{0.9\textwidth}
\textbf{Algorithm:}
\begin{itemize}
\item[\bf 1.]
Let $(\widehat X_1, \widehat X_2) = 0.5(I, I)$ be initial matrix solution of iteration and $\varepsilon = 0.001$.

\item[\bf 2.]
Substituting $(X_1, X_2) = (\widehat X_1, \widehat X_2)$ into (\ref{eq:cond}) yields $(\Theta_1, \Theta_2) = (\Theta_1^*, \Theta_2^*)$. 

\item[\bf 3.]
Set 
\begin{equation}\label{cA-1}
\cA_1\deq A+B_2\Theta_2^*,  \quad  
\cQ_1\deq Q_1+(\Theta_2^*)^\top R_{122}\Theta_2^*, \quad 
\cS_{11}\deq R_{112}\Theta_2^*. 
\end{equation}
Then the first equation in (\ref{eq:Riccati}) and the first equation in (\ref{eq:cond}) can be reduced to 
\begin{equation}\label{ARE-1}
X_1\cA_1+\cA_1^\top X_1+\cQ_1-\big(X_1B_1+\cS_{11}^\top\big)R_{111}^{-1}\big(B_1^\top X_1+\cS_{11}\big)=0,\\
\end{equation}
and 
\begin{equation}\label{Th1-Th*1-closed}
B_1^\top X_1+\cS_{11}+R_{111}\Theta_1=0,
\end{equation}
respectively. Solving (\ref{ARE-1}) yields its solution $X_1^*$, and substituting $X_1 = X_1^*$ into (\ref{Th1-Th*1-closed}) leads to $\Theta_1^*$. 

%
\end{itemize}
\end{boxedminipage}

\ms

\begin{boxedminipage}[t]{0.9\textwidth}
\textbf{Algorithm:}
\begin{itemize}
%
%
%
\item[\bf 4.] 
Similarly, set 
\begin{equation}\label{cA-2}
\cA_2\deq A+B_1\Theta_1^*,  \quad  
\cQ_2\deq Q_2+(\Theta_1^*)^\top R_{211}\Theta_1^*, \quad 
\cS_{22}\deq R_{221}\Theta_1^*. 
\end{equation}
Then the second equation in (\ref{eq:Riccati}) and the second equation in (\ref{eq:cond}) can be reduced to 
\begin{equation}\label{ARE-2}
X_2\cA_2+\cA_2^\top X_2+\cQ_2-\big(X_2B_2+\cS_{22}^\top\big)R_{222}^{-1}\big(B_2^\top X_2+\cS_{22}\big)=0,\\
\end{equation}
and 
\begin{equation}\label{Th1-Th*2-closed}
B_2^\top X_2+\cS_{22}+R_{222}\Theta_2=0,  
\end{equation}
respectively. Solving (\ref{ARE-2}) yields its solution $X_2^*$, and substituting $X_2 = X_2^*$ into (\ref{Th1-Th*2-closed}) leads to $\Theta_2^*$. 

\item[\bf 5.]
If $\Vert\widehat X_1 - X_1^*\Vert^2 + \Vert\widehat X_2 - X_2^*\Vert^2 \ges \varepsilon$, 
then set $(\widehat X_1, \widehat X_2 )= (X_1^*, X_2^*)$ and go to Step \textbf{2}. Otherwise, stop.
\end{itemize}
\end{boxedminipage}

\ms

\ms

Since $Q_i \succeq 0$ and $R_i \succ 0$ for $i=1,2$, we have $\cQ_i \succeq 0$. In addition, $R_{111} \succeq 0$ and $R_{222} \succeq 0$. Therefore, applying the Schur method in \cite{Laub}, we can tackle equations (\ref{ARE-1}) and (\ref{ARE-2}) directly. Although the Schur-based method has all the reliability and efficiency, it is relatively difficult to effectively parallelize and vectorize. Therefore, other methods have been re-examined in order to implement them on various advanced computing architectures.
Due to its non-linearity in $(X_1, X_2)$, it is difficult to study its numerical analysis. Numerical theory will be the subject of our future research.
In addition, we also can tackle equations (\ref{ARE-1}) and (\ref{ARE-2}) using 
a computational approach via a semi-definite programming associated with linear matrix inequalities, whose feasibility is  equivalent to the solvability of equations (\ref{ARE-1}) and (\ref{ARE-2}).

\section{Examples}

In this section, we present some examples illustrating the algorithm in the previous section. 
We discuss about five cases: 
(1) $n=1$, $m_1=1$ and $m_2=1$; 
(2) $n=2$, $m_1=1$ and $m_2=1$;
(3) $n=2$, $m_1=2$ and $m_2=1$;
(4) $n=2$, $m_1=2$ and $m_2=2$;
(5) $n=3$, $m_1=2$ and $m_2=2$.
All computations are done on an iMac using Matlab software and double precision arithmetic.

\ms

\bex{Example 1} \rm 
Consider the following problem with $n=1$, $m_1=1$ and $m_2=1$.
In this example,
$$\left\{\begin{array}{ll}
A=1, \quad B_1=B_2=1, \quad Q_1=2, \quad Q_2=2.5, \\ [2mm] 
R_1=\begin{pmatrix} 1 & 0 \\ 0 & 0\end{pmatrix}, \quad R_2=\begin{pmatrix} 0 & 0 \\ 0 & 2\end{pmatrix}.
\end{array}\right.$$
Then the corresponding system of generalized AREs (\ref{eq:Riccati})-(\ref{eq:cond}) reads
\begin{equation}\label{ex1:Riccati}
\left\{\begin{array}{ll}
X_1+X_1+\Theta^\top R_1\Theta + X_1(1,1)\Theta + \Theta^\top (1,1)^\top X_1 + 2 = 0, & X_1+\Theta_1=0, \\ [2mm]
X_2+X_2+\Theta^\top R_2\Theta + X_2(1,1)\Theta + \Theta^\top (1,1)^\top X_2 + 2.5 = 0, & X_2+2\Theta_2=0,
\end{array}\right.
\end{equation}
which is equivalent to 
\begin{equation}\label{ex1:Riccati1}
\left\{\begin{array}{ll}
2X_1 - X_1^2 - X_1X_2 + 2 = 0, \\ [2mm]
2X_2 - \frac{1}{2}X_2^2 - 2X_1X_2 + 2.5 = 0.
\end{array}\right.
\end{equation}
This is an example with two scalar variables and two quadratic equations.
Let $(\widehat X_1, \widehat X_2) = (0.5, 0.5)$ be initial solution of iteration and $\varepsilon = 0.001$. Using
algorithm in Section \ref{sect:algorithm}, we can get 
\begin{equation*}
\left\{\begin{array}{lll} 
\cA_1 = 1+\Theta_2^*,  \quad  \cQ_1 = 2, \quad \cS_{11} = 0, \\ [2mm]
\cA_2 = 1+\Theta_1^*,  \quad  \cQ_2 = 2.5, \quad \cS_{22} = 0.
\end{array}\right.
\end{equation*}
Then (\ref{ARE-1})-(\ref{Th1-Th*1-closed}) and (\ref{ARE-2})-(\ref{Th1-Th*2-closed}) reads
\begin{equation*}
\left\{\begin{array}{lll} 
X_1\cA_1+\cA_1X_1+2-X_1^2=0, \quad X_1+\Theta_1=0, \\ [2mm]
X_2\cA_2+\cA_2X_2+2.5-\frac{1}{2}X_2^2=0, \quad X_2+2\Theta_2=0.
\end{array}\right.
\end{equation*}
Several iterations leads to a pair of solution $(X_1^*, X_2^*) = (2, 1)$ satisfying (\ref{ex1:Riccati}) and its equivalent equations (\ref{ex1:Riccati1}).
\ex


\ms

\bex{Example 2} \rm 
Consider the following problem with $n=2$, $m_1=1$ and $m_2=1$.
In this example,
$$\left\{\begin{array}{ll}
A=\begin{pmatrix} -1 & 0 \\ 0 & -1\end{pmatrix}, \quad 
B = (B_1, B_2), \quad B_1=\begin{pmatrix} 1 \\ 0 \end{pmatrix}, \quad B_2=\begin{pmatrix} 0 \\ 1 \end{pmatrix}, \\ [5mm] 
Q_1=\begin{pmatrix} 3 & 0 \\ 0 & 2.75\end{pmatrix}, \quad Q_2=\begin{pmatrix} 2 & 0 \\ 0 & 2.5\end{pmatrix}, \quad 
R_1=\begin{pmatrix} 1 & 0 \\ 0 & 1\end{pmatrix}, \quad R_2=\begin{pmatrix} 2 & 0 \\ 0 & 2\end{pmatrix}.
\end{array}\right.$$
Then the corresponding system of generalized AREs (\ref{eq:Riccati})-(\ref{eq:cond}) reads
\begin{equation}\label{ex2:Riccati}
\left\{\begin{array}{ll}
X_1A+A^\top X_1+\Theta^\top R_1\Theta + X_1B\Theta + \Theta^\top B^\top X_1 + Q_1 = 0, \\ [2mm]
X_2A+A^\top X_2+\Theta^\top R_2\Theta + X_2B\Theta + \Theta^\top B^\top X_2 + Q_2 = 0, \\ [2mm]
\begin{pmatrix} B_1^\top X_1 \\ B_2^\top X_2\end{pmatrix}
+ \begin{pmatrix} 1 & 0 \\ 0 & 2\end{pmatrix} 
\begin{pmatrix}\Theta_1 \\ \Theta_2\end{pmatrix} = 0, 
\end{array}\right.
\end{equation}
where $\Theta = (\Theta_1, \Theta_2)^\top$, $\Theta_1 \in \mathbb{R}^{1\times 2}$ and $\Theta_2 \in \mathbb{R}^{1\times 2}$.
Let $(\widehat X_1, \widehat X_2)$ be initial matrix solution of iteration, where
\begin{equation}\label{ex2:initial}
\widehat X_1 = \begin{pmatrix} 0.5 & 0 \\ 0 & 0.5\end{pmatrix}, \quad 
\widehat X_2 = \begin{pmatrix} 0.5 & 0 \\ 0 & 0.5\end{pmatrix},
\end{equation}
and $\varepsilon = 0.001$. Using
algorithm in Section \ref{sect:algorithm}, we can get 
\begin{equation*}
\left\{\begin{array}{lll} 
\cA_1 = A+B_2\Theta_2^*,  \quad  \cQ_1 = Q_1, \quad \cS_{11} = 0, \\ [2mm]
\cA_2 = A+B_1\Theta_1^*,  \quad  \cQ_2 = Q_2, \quad \cS_{22} = 0.
\end{array}\right.
\end{equation*}
Then (\ref{ARE-1})-(\ref{Th1-Th*1-closed}) and (\ref{ARE-2})-(\ref{Th1-Th*2-closed}) reads 
\begin{equation*}
\left\{\begin{array}{lll} 
X_1\cA_1+\cA_1^\top X_1+Q_1-X_1B_1(B_1^\top X_1)=0, \\ [2mm]
X_2\cA_2+\cA_2^\top X_2+Q_2-\frac{1}{2}X_2B_2(B_2^\top X_2)=0, \\ [2mm]
\begin{pmatrix} B_1^\top X_1 \\ B_2^\top X_2\end{pmatrix}
+ \begin{pmatrix} 1 & 0 \\ 0 & 2\end{pmatrix} 
\begin{pmatrix}\Theta_1 \\ \Theta_2\end{pmatrix} = 0.
\end{array}\right.
\end{equation*}
Several iterations leads to a pair of solution 
\begin{equation}\label{ex2:sol}
X_1^* = \begin{pmatrix} 1 & 0 \\ 0 & 1\end{pmatrix} \quad \mbox{and} \quad
X_2^* = \begin{pmatrix} 1 & 0 \\ 0 & 1\end{pmatrix}
\end{equation}
satisfying (\ref{ex2:Riccati}).

\ex


\ms

\bex{Example 3} \rm 
Consider the following problem with $n=2$, $m_1=2$ and $m_2=1$.
In this example,
$$\left\{\begin{array}{ll}
A=\begin{pmatrix} -1 & 0 \\ 0 & -1\end{pmatrix}, \quad 
B = (B_1, B_2), \quad B_1=\begin{pmatrix} 1 & 0 \\ 0 & 1\end{pmatrix}, \quad B_2=\begin{pmatrix} 1 \\ 1 \end{pmatrix}, \\ [5mm] 
Q_1=\begin{pmatrix} 4.25 & 2.25 \\ 2.25 & 3\end{pmatrix}, \quad Q_2=\begin{pmatrix} 2.5 & 0 \\ 0 & 2\end{pmatrix}, \quad 
R_1=\begin{pmatrix} 1 & 0.5 & 0 \\ 0.5 & 1 & 0 \\ 0 & 0 & 1\end{pmatrix}, \quad R_2=\begin{pmatrix} 2 & 0 & 0 \\ 0 & 2 & 1 \\ 0 & 1 & 2\end{pmatrix}.
\end{array}\right.$$
Then the corresponding system of generalized AREs (\ref{eq:Riccati})-(\ref{eq:cond}) reads
\begin{equation}\label{ex3:Riccati}
\left\{\begin{array}{ll}
X_1A+A^\top X_1+\Theta^\top R_1\Theta + X_1B\Theta + \Theta^\top B^\top X_1 + Q_1 = 0, \\ [2mm]
X_2A+A^\top X_2+\Theta^\top R_2\Theta + X_2B\Theta + \Theta^\top B^\top X_2 + Q_2 = 0, \\ [2mm]
\begin{pmatrix} B_1^\top X_1 \\ B_2^\top X_2\end{pmatrix}
+ \begin{pmatrix} 1 & 0.5 & 0 \\ 0.5 & 1 & 0 \\ 0 & 1 & 2\end{pmatrix}
\begin{pmatrix}\Theta_1 \\ \Theta_2\end{pmatrix} = 0, 
\end{array}\right.
\end{equation}
where $\Theta = (\Theta_1, \Theta_2)^\top$, $\Theta_1 \in \mathbb{R}^{2\times 2}$ and $\Theta_2 \in \mathbb{R}^{1\times 2}$.
Let $(\widehat X_1, \widehat X_2)$ be initial matrix solution of iteration, where
\begin{equation}\label{ex3:initial}
\widehat X_1 = \begin{pmatrix} 0.5 & 0 \\ 0 & 0.5\end{pmatrix}, \quad 
\widehat X_2 = \begin{pmatrix} 0.5 & 0 \\ 0 & 0.5\end{pmatrix},
\end{equation}
and $\varepsilon = 0.001$. Using
algorithm in Section \ref{sect:algorithm}, we can get 
\begin{equation*}
\left\{\begin{array}{lll} 
\cA_1 = A+B_2\Theta_2^*,  &  \cQ_1 = Q_1+(\Theta_2^*)^\top\Theta_2^*, & 
\cS_{11} = 0, \\ [2mm]
\cA_2 = A+B_1\Theta_1^*,  &  \cQ_2 = Q_2+(\Theta_1^*)^\top R_{211}\Theta_1^*, & 
\cS_{22} = R_{221}\Theta_1^*,
\end{array}\right.
\end{equation*}
where $R_{211} = \begin{pmatrix} 2 & 0 \\ 0 & 2\end{pmatrix}$ and $R_{221} = \begin{pmatrix} 0 & 1\end{pmatrix}$. 
Then (\ref{ARE-1})-(\ref{Th1-Th*1-closed}) and (\ref{ARE-2})-(\ref{Th1-Th*2-closed}) reads 
\begin{equation*}
\left\{\begin{array}{lll} 
X_1\cA_1+\cA_1^\top X_1+\cQ_1-X_1B_1R_{111}^{-1}(B_1^\top X_1)=0, \\ [2mm]
X_2\cA_2+\cA_2^\top X_2+\cQ_2-\frac{1}{2}(X_2B_2+\cS_{22}^\top)(B_2^\top X_2+\cS_{22})=0, \\ [2mm]
\begin{pmatrix} B_1^\top X_1 \\ B_2^\top X_2\end{pmatrix}
+ \begin{pmatrix} 1 & 0.5 & 0 \\ 0.5 & 1 & 0 \\ 0 & 1 & 2\end{pmatrix}
\begin{pmatrix}\Theta_1 \\ \Theta_2\end{pmatrix} = 0,
\end{array}\right.
\end{equation*}
where $R_{111} = \begin{pmatrix} 1 & 0.5 \\ 0.5 & 1\end{pmatrix}$.
Several iterations leads to a pair of solution 
\begin{equation}\label{ex3:sol}
X_1^* = \begin{pmatrix} 1 & 0.5 \\ 0.5 & 1\end{pmatrix} \quad \mbox{and} \quad
X_2^* = \begin{pmatrix} 1 & 0 \\ 0 & 1\end{pmatrix}
\end{equation}
satisfying (\ref{ex3:Riccati}).

\ex


\ms

\bex{Example 4} \rm 
Consider the following problem with $n=2$, $m_1=2$ and $m_2=2$.
In this example,
$$\left\{\begin{array}{ll}
A=\begin{pmatrix} -1 & 0 \\ 0 & -1\end{pmatrix}, \quad 
B = (B_1, B_2), \quad B_1=\begin{pmatrix} 1 & 0 \\ 0 & 1\end{pmatrix}, \quad B_2=\begin{pmatrix} 1 & 0.5 \\ 0.5 & 1 \end{pmatrix}, \\ [5mm] 
Q_1=\begin{pmatrix} 4\frac{1}{2} & 2\frac{3}{4} \\ 2\frac{3}{4} & 3\frac{1}{6}\end{pmatrix}, \quad Q_2=\begin{pmatrix} 8 & -\frac{1}{2} \\ -\frac{1}{2} & 3\frac{1}{6}\end{pmatrix}, \\ [5mm]
R_1=\begin{pmatrix} 1 & 0.5 & 0 & 0 \\ 0.5 & 1 & 0 & 0 \\ 0 & 0 & 1 & 0.5 \\ 0 & 0 & 0.5 & 1\end{pmatrix}, \quad R_2=\begin{pmatrix} 2 & 0 & 0 & 0 \\ 0 & 2 & 1 & 0 \\ 0 & 1 & 2 & 1 \\ 0 & 0 & 1 & 2\end{pmatrix}.
\end{array}\right.$$
Then the corresponding system of generalized AREs (\ref{eq:Riccati})-(\ref{eq:cond}) reads
\begin{equation}\label{ex4:Riccati}
\left\{\begin{array}{ll}
X_1A+A^\top X_1+\Theta^\top R_1\Theta + X_1B\Theta + \Theta^\top B^\top X_1 + Q_1 = 0, \\ [2mm]
X_2A+A^\top X_2+\Theta^\top R_2\Theta + X_2B\Theta + \Theta^\top B^\top X_2 + Q_2 = 0, \\ [2mm]
\begin{pmatrix} B_1^\top X_1 \\ B_2^\top X_2\end{pmatrix}
+ \begin{pmatrix} 1 & 0.5 & 0 & 0 \\ 0.5 & 1 & 0 & 0 \\ 0 & 1 & 2 & 1 \\ 0 & 0 & 1 & 2\end{pmatrix}
\begin{pmatrix}\Theta_1 \\ \Theta_2\end{pmatrix} = 0, 
\end{array}\right.
\end{equation}
where $\Theta = (\Theta_1, \Theta_2)^\top$, $\Theta_1 \in \mathbb{R}^{2\times 2}$ and $\Theta_2 \in \mathbb{R}^{2\times 2}$.
Let $(\widehat X_1, \widehat X_2)$ be initial matrix solution of iteration, where
\begin{equation}\label{ex4:initial}
\widehat X_1 = \begin{pmatrix} 0.5 & 0 \\ 0 & 0.5\end{pmatrix}, \quad 
\widehat X_2 = \begin{pmatrix} 0.5 & 0 \\ 0 & 0.5\end{pmatrix},
\end{equation}
and $\varepsilon = 0.001$. Using
algorithm in Section \ref{sect:algorithm}, we can get 
\begin{equation*}
\left\{\begin{array}{lll} 
\cA_1 = A+B_2\Theta_2^*,  &  \cQ_1 = Q_1+(\Theta_2^*)^\top R_{122}\Theta_2^*, & 
\cS_{11} = 0, \\ [2mm]
\cA_2 = A+B_1\Theta_1^*,  &  \cQ_2 = Q_2+(\Theta_1^*)^\top R_{211}\Theta_1^*, & 
\cS_{22} = R_{221}\Theta_1^*,
\end{array}\right.
\end{equation*}
where $R_{122} = \begin{pmatrix} 1 & 0.5 \\ 0.5 & 1\end{pmatrix}$, $R_{211} = \begin{pmatrix} 2 & 0 \\ 0 & 2\end{pmatrix}$ and $R_{221} = \begin{pmatrix} 0 & 1 \\ 0 & 0\end{pmatrix}$. 
Then (\ref{ARE-1})-(\ref{Th1-Th*1-closed}) and (\ref{ARE-2})-(\ref{Th1-Th*2-closed}) reads 
\begin{equation*}
\left\{\begin{array}{lll} 
X_1\cA_1+\cA_1^\top X_1+\cQ_1-X_1B_1R_{111}^{-1}(B_1^\top X_1)=0, \\ [2mm]
X_2\cA_2+\cA_2^\top X_2+\cQ_2-(X_2B_2+\cS_{22}^\top)R_{222}^{-1}(B_2^\top X_2+\cS_{22})=0, \\ [2mm]
\begin{pmatrix} B_1^\top X_1 \\ B_2^\top X_2\end{pmatrix}
+ \begin{pmatrix} 1 & 0.5 & 0 & 0 \\ 0.5 & 1 & 0 & 0 \\ 0 & 1 & 2 & 1 \\ 0 & 0 & 1 & 2\end{pmatrix}
\begin{pmatrix}\Theta_1 \\ \Theta_2\end{pmatrix} = 0,
\end{array}\right.
\end{equation*}
where $R_{111} = \begin{pmatrix} 1 & 0.5 \\ 0.5 & 1\end{pmatrix}$ and $R_{222} = \begin{pmatrix} 2 & 1 \\ 1 & 2\end{pmatrix}$.
Several iterations leads to a pair of solution 
\begin{equation}\label{ex4:sol}
X_1^* = \begin{pmatrix} 1 & 0.5 \\ 0.5 & 1\end{pmatrix} \quad \mbox{and} \quad 
X_2^* = \begin{pmatrix} 2 & 0 \\ 0 & 1\end{pmatrix}
\end{equation}
satisfying (\ref{ex4:Riccati}).

\ex


\ms

\bex{Example 5} \rm 
Consider the following problem with $n=3$, $m_1=2$ and $m_2=2$.
In this example,
$$\left\{\begin{array}{ll}
A=\begin{pmatrix} -5 & 0 & -1 \\ 0 & -10 & 0 \\ -1 & 0 & -5\end{pmatrix}, \quad 
B = (B_1, B_2), \quad B_1=\begin{pmatrix} 1 & 0 \\ 0 & 1 \\ 0 & 1\end{pmatrix}, \quad B_2=\begin{pmatrix} 1 & 0.5 \\ 0.5 & 1 \\ 0 & 1\end{pmatrix}, \\ [8mm] 
Q_1=\begin{pmatrix} 13.1852 &  19.1111  &  3.8704 \\
   19.1111 & 43.4167  &  3.9722 \\
    3.8704 &  3.9722  & 11.8241\end{pmatrix}, \quad 
Q_2=\begin{pmatrix} 21.8519  & -0.7222  &  3.9815 \\
   -0.7222 & 26.1667  & 25.2222 \\
    3.9815  & 25.2222  & 31.6852\end{pmatrix}, \\ [8mm]
R_1=\begin{pmatrix} 1 & 0.5 & 0 & 0 \\ 0.5 & 1 & 0 & 0 \\ 0 & 0 & 1 & 0.5 \\ 0 & 0 & 0.5 & 1\end{pmatrix}, \quad R_2=\begin{pmatrix} 2 & 0 & 0 & 0 \\ 0 & 2 & 1 & 0 \\ 0 & 1 & 2 & 1 \\ 0 & 0 & 1 & 2\end{pmatrix}.
\end{array}\right.$$
Then the corresponding system of generalized AREs (\ref{eq:Riccati})-(\ref{eq:cond}) reads
\begin{equation}\label{ex5:Riccati}
\left\{\begin{array}{ll}
X_1A+A^\top X_1+\Theta^\top R_1\Theta + X_1B\Theta + \Theta^\top B^\top X_1 + Q_1 = 0, \\ [2mm]
X_2A+A^\top X_2+\Theta^\top R_2\Theta + X_2B\Theta + \Theta^\top B^\top X_2 + Q_2 = 0, \\ [2mm]
\begin{pmatrix} B_1^\top X_1 \\ B_2^\top X_2\end{pmatrix}
+ \begin{pmatrix} 1 & 0.5 & 0 & 0 \\ 0.5 & 1 & 0 & 0 \\ 0 & 1 & 2 & 1 \\ 0 & 0 & 1 & 2\end{pmatrix}
\begin{pmatrix}\Theta_1 \\ \Theta_2\end{pmatrix} = 0, 
\end{array}\right.
\end{equation}
where $\Theta = (\Theta_1, \Theta_2)^\top$, $\Theta_1 \in \mathbb{R}^{2\times 3}$ and $\Theta_2 \in \mathbb{R}^{2\times 3}$.
Let $(\widehat X_1, \widehat X_2)$ be initial matrix solution of iteration, where
\begin{equation}\label{ex5:initial}
\widehat X_1 = \begin{pmatrix} 0.5 & 0 & 0 \\ 0 & 0.5 & 0 \\ 0 & 0 & 0.5\end{pmatrix}, \quad 
\widehat X_2 = \begin{pmatrix} 0.5 & 0 & 0 \\ 0 & 0.5 & 0 \\ 0 & 0 & 0.5\end{pmatrix},
\end{equation}
and $\varepsilon = 0.001$. Using
algorithm in Section \ref{sect:algorithm}, we can get 
\begin{equation*}
\left\{\begin{array}{lll} 
\cA_1 = A+B_2\Theta_2^*,  &  \cQ_1 = Q_1+(\Theta_2^*)^\top R_{122}\Theta_2^*, & 
\cS_{11} = 0, \\ [2mm]
\cA_2 = A+B_1\Theta_1^*,  &  \cQ_2 = Q_2+(\Theta_1^*)^\top R_{211}\Theta_1^*, & 
\cS_{22} = R_{221}\Theta_1^*,
\end{array}\right.
\end{equation*}
where $R_{122} = \begin{pmatrix} 1 & 0.5 \\ 0.5 & 1\end{pmatrix}$, $R_{211} = \begin{pmatrix} 2 & 0 \\ 0 & 2\end{pmatrix}$ and $R_{221} = \begin{pmatrix} 0 & 1 \\ 0 & 0\end{pmatrix}$. 
Then (\ref{ARE-1})-(\ref{Th1-Th*1-closed}) and (\ref{ARE-2})-(\ref{Th1-Th*2-closed}) reads 
\begin{equation*}
\left\{\begin{array}{lll} 
X_1\cA_1+\cA_1^\top X_1+\cQ_1-X_1B_1R_{111}^{-1}(B_1^\top X_1)=0, \\ [2mm]
X_2\cA_2+\cA_2^\top X_2+\cQ_2-(X_2B_2+\cS_{22}^\top)R_{222}^{-1}(B_2^\top X_2+\cS_{22})=0, \\ [2mm]
\begin{pmatrix} B_1^\top X_1 \\ B_2^\top X_2\end{pmatrix}
+ \begin{pmatrix} 1 & 0.5 & 0 & 0 \\ 0.5 & 1 & 0 & 0 \\ 0 & 1 & 2 & 1 \\ 0 & 0 & 1 & 2\end{pmatrix}
\begin{pmatrix}\Theta_1 \\ \Theta_2\end{pmatrix} = 0,
\end{array}\right.
\end{equation*}
where $R_{111} = \begin{pmatrix} 1 & 0.5 \\ 0.5 & 1\end{pmatrix}$ and $R_{222} = \begin{pmatrix} 2 & 1 \\ 1 & 2\end{pmatrix}$.
Several iterations leads to a pair of solution 
\begin{equation}\label{ex5:sol}
X_1^* = \begin{pmatrix} 1 & 1 & 0 \\ 1 & 2 & 0 \\ 0 & 0 & 1\end{pmatrix} \quad \mbox{and} \quad 
X_2^* = \begin{pmatrix} 2 & 0 & 0 \\ 0 & 1 & 1 \\ 0 & 1 & 2\end{pmatrix}
\end{equation}
satisfying (\ref{ex5:Riccati}).

\ex

\section{Conclusion}

We have learned many methods of algebraic Riccati equations and their performance in computer finite arithmetic environments. However, many important questions remain unanswered and are the subject of ongoing research topics.


\ms

\end{document}